\documentclass{article}
\usepackage{amsmath,amssymb,bm}

\title{Linear projections and successive minima}
\author{Christophe Soul\'e}
\date{}

\begin{document}

\maketitle

Let $K$ be a number field, ${\mathcal O}_K$ its ring of integers and $E$ a projective ${\mathcal O}_K$-module of finite rank $N$. We endow $E \underset{\mathbb Z}{\otimes} {\mathbb C}$ with an hermitian metric $h$ and we let $\mu_1 , \ldots , \mu_N$ be the logarithm of the successive minima of $(E,h)$. Assume $X_K \subset {\mathbb P} (E_K^{\vee})$ is a smooth geometrically irreducible curve. In this paper we shall find a lower bound for the numbers $\mu_i$, $3 \leq i \leq N$, in terms of the height of $X_K$, $\mu_1$ and the average of the $\mu_i$'s (Theorem~2). This result is a complement to \cite{S}, Theorem~4, which gives a lower bound for $\mu_1$.

\smallskip

The method of proof is a variant of \cite{S}, {\it loc. cit.} It relies upon Morrison's proof of the fact that $X_K$ is Chow semi-stable \cite{Mo}. We use a filtration $V_1 = E_K \supset V_2 \supset \ldots \supset V_N$ of the vector space $E_K$. This filtration is chosen so that, for suitable values of $i$, the projection ${\mathbb P} (V_i^{\vee}) \cdots \to {\mathbb P} (V_{i+1}^{\vee})$ does not change the degree of the image of $X_K$ by linear projection. That such a choice is possible follows from a result of C.~Voisin, namely an effective version of a theorem of Segre on linear projections of complex projective curves (Theorem~1). I thank her for proving this result and for helpful discussions.

\section{Linear projections of projective curves}

Let $C \subset {\mathbb P}^n$ be an integral projective curve over ${\mathbb C}$ and $d$ its degree. Assume that $C$ is not contained in some hyperplane, $d \geq 3$ and $n \geq 3$.

\bigskip

\noindent {\bf Theorem 1.} (C. Voisin) {\it There exists an integer $A = A(d)$ and a finite set $\Sigma$ of points in $({\mathbb P}^n - C) ({\mathbb C})$, of order at most $A$, such that, for every point $P \in {\mathbb P}^n ({\mathbb C}) - \Sigma \cup C({\mathbb C})$, the linear projection ${\mathbb P}^n \cdots \to {\mathbb P}^{n-1}$ of center $P$ maps $C$ birationally onto its image.}

\bigskip

\noindent {\it Proof.} The existence of a finite set $\Sigma$ with the property above is a special case of a theorem of C.~Segre \cite{C}. The order of $\Sigma$ can be bounded as follows by a function of $d$.

\smallskip

If $n > 3$ a generic linear projection into ${\mathbb P}^3$ will map $C$ isomorphically onto its image \cite{H} and the exceptional set $\Sigma \subset {\mathbb P}^n$ bijectively onto the exceptional set in ${\mathbb P}^3$. Therefore we can assume that $n=3$.

\smallskip

When the projection with center $P \in {\mathbb P}^3 ({\mathbb C})$ is not birational from the curve $C$ to its image $C' \subset {\mathbb P}^2$, we have $d' = \deg (C') \leq \frac{d}{2}$ hence $d' \leq d-2$, and $P$ is the vertex of a cone $K$ with base $C'$ containing $C$. So we have to bound the number of such cones.

\smallskip

Let $N$ be the dimension of the kernel of the restriction map
$$
\alpha : H^0 ({\mathbb P}^3 , {\mathcal O} (d')) \to H^0 (C , {\mathcal O} (d')) \, .
$$
Clearly $N$ is bounded as a function of $d$ and any $f \in \ker (\alpha)$ is an homogeneous polynomial of degree $d'$ which vanishes on $C$.

\smallskip

Let $Z \subset {\mathbb P}^3 ({\mathbb C}) \times {\mathbb P}^{N-1} ({\mathbb C})$ be the set of pairs $(P,f)$ such that $f$ is the equation of a cone $K$ of vertex $P$. If $p_1 : {\mathbb P}^3 \times {\mathbb P}^{N-1} \to {\mathbb P}^3$ is the first projection, we have to bound the order of $p_1 (Z)$. We note that this order is at most the number $c$ of connected components of $Z$.

\smallskip

Now $Z$ is defined by equations of bidegree $(\delta , 1)$, $\delta \leq d'$. Indeed $f$ is homogeneous of degree $d'$ and $(P,f) \in Z$ when all the derivatives of $f$, except those of order $d'$, vanish at $P$.

\smallskip

Let $L = {\mathcal O} (d',1)$, $M = \dim H^0 ({\mathbb P}^3 \times {\mathbb P}^N , L) - 1$, and
$$
j : {\mathbb P}^3 \times {\mathbb P}^N \to {\mathbb P}^M
$$
the Segre embedding. Since $j(Z)$ is the intersection of $j({\mathbb P}^3 \times {\mathbb P}^N)$ with linear hyperplanes, B\'ezout theorem (\cite{F}, \S~8.4) tells us that
$$
c \leq \deg (j({\mathbb P}^3 \times {\mathbb P}^N)) \, .
$$
Hence $c$ is bounded by a function of $d$. \hfill $\Box$

\bigskip

\noindent {\bf Corollary.} {\it Given any projective line $\Lambda \subset {\mathbb P}^n$, there exists a finite set $\Phi$ of order at most $A(d) + d$ in $\Lambda$ such that, if $P \in \Lambda - \Phi$, the linear projection of center $P$ maps $C$ birationally onto its image.}

\bigskip

\noindent {\it Proof.} Since $C$ is not equal to $\Lambda$, the cardinality of $C \cap \Lambda$ is at most $d$. So the Corollary follows from Theorem~1.

\section{Successive minima}

\subsection{ \ }

Let $K$ be a number field, $[K:{\mathbb Q}]$ its degree over ${\mathbb Q}$, ${\mathcal O}_K$ its ring of integers, $S = {\rm Spec} ({\mathcal O}_K)$ the associated scheme and $\Sigma$ the set of complex embeddings of $K$. Consider an hermitian vector bundle $(E,h)$ over $S$, {\it i.e.} $E$ is a torsion free ${\mathcal O}_K$-module of finite rank $N$ and, for all $\sigma \in \Sigma$, the associated complex vector space $E_{\sigma} = E \underset{{\mathcal O}_K}{\otimes} {\mathbb C}$ is equipped with an hermitian scalar product $h_{\sigma}$. If $\bar\sigma$ is the conjugate of $\sigma$, we assume that the complex conjugation $E_{\sigma} \simeq E_{\bar\sigma}$ is an isometry.

\smallskip

If $i$ is a positive integer, $i \leq N$, we let $\mu_i$ be the infimum of the set of real numbers $r$ such that there exist $v_1 , \ldots , v_i \in E$, linearly independent over $K$, such that $\log \Vert v_{\alpha} \Vert \leq r$ for all $\alpha \leq i$. The number $\mu_i$ is thus the logarithm of the $i$-th successive minimum of $(E,h)$. Let
\begin{equation}
\label{eq1}
\mu = \frac{\mu_1 + \cdots + \mu_N}{N} \, .
\end{equation}

\subsection{}

If $E^{\vee} = {\rm Hom} (E , {\mathcal O}_K)$ is the dual of $E$ we let ${\mathbb P} (E^{\vee})$ be the associated projective space, representing lines in $E^{\vee}$. Let $E_K^{\vee} = E^{\vee} \underset{{\mathcal O}_K}{\otimes} K$ and $X_K \subset {\mathbb P} (E_K^{\vee})$ a smooth geometrically irreducible curve of genus $g$ and degree $d$. We assume that the embedding of $X_K$ into ${\mathbb P} (E_K^{\vee})$ is defined by a complete linear series on $X_K$. We also assume that $g \geq 2$ and $d \geq 2g+1$. The rank  of $E$ is thus $N = d+1-g$.

\smallskip

If $X$ is the Zariski closure of $X_K$ in ${\mathbb P} (E^{\vee})$ and $\overline{{\mathcal O} (1)}$ the canonical hermitian line bundle on ${\mathbb P} (E^{\vee})$, the Faltings height of $X_K$ is the real number
$$
h (X_K) = \widehat\deg \, (\hat c_1 (\overline{{\mathcal O} (1)})^2 \mid X) \, ,
$$
see \cite{BGS} (3.1.1) and (3.1.5).

\subsection{}

For any positive integer $i \leq N$ we define the integer $f_i$ by the formulae
$$
f_i = i-1 \qquad \mbox{if} \qquad i-1 \leq d - 2g
$$
and
$$
f_i = i-1 + \alpha \qquad \mbox{if} \qquad i-1 = d -2g + \alpha \, , \qquad 0 \leq \alpha \leq g \, .
$$

Assume $k$ and $i$ are two positive integers, $k \leq N$, $i \leq N$. We let
$$
h_{i,k} = \left\{ \begin{matrix}
f_i &\mbox{if} &i \leq k \, , \ i = N-1 \ \mbox{or} \ i = N \\
f_k &\mbox{if} &k \leq i \leq N-2 \, . \hfill
\end{matrix} \right.
$$
Finally, if $2 \leq k \leq N$, we let
$$
B_k = \underset{i=2,\ldots , N}{\rm max} \ \frac{h_{i,k}^2}{(i-1) \, h_{i,k} - \underset{j=1}{\overset{i-1}{\sum}} \, h_{j,k}} \, .
$$

\bigskip

\noindent {\bf Theorem 2.} {\it There exists a constant $C = C(d)$ such that, for every $k$ such that $2 \leq k \leq N-3$} ,
$$
B_k (\mu_{N+1-k} - \mu) + \frac{h(X_K)}{[K:{\mathbb Q}]} + 2d \, \mu \geq (2d - (N+1) \, B_k) (\mu - \mu_1) - C \, .
$$

\subsection{}

To prove Theorem 2 fix a positive integer $k \leq N-3$ and choose elements $x_1 , \ldots , x_N$ in $E$, linearly independent over $K$ and such that
$$
\log \Vert x_i \Vert = \mu_{N-i+1} \, , \quad 1 \leq i \leq N \, .
$$

Fix integers $n_{\alpha}$, $\alpha = k+1 , \ldots , N-2$, to be specified later (in \S~2.6). If $1 \leq i \leq N$ we define
\begin{equation}
\label{eq2}
v_i = \left\{ \begin{matrix}
x_i + n_i \, x_{i-1} &\mbox{if} \ k+1 \leq i \leq N-2 \\
x_i \hfill &\mbox{else.} \hfill
\end{matrix} \right.
\end{equation}
We get a complete flag $E_K = V_1 \supset V_2 \supset \ldots \supset V_N$ by defining $V_i$ to be the linear span of $v_i , v_{i+1} , \ldots , v_N$.

\smallskip

When $m$ is large enough the cup-product map
$$
\varphi : E_K^{\otimes m} \to H^0 (X_K , {\mathcal O} (m))
$$
is surjective, hence $H^0 (X_K , {\mathcal O} (m))$ is generated by the monomials
$$
v_1^{\alpha_1} \ldots v_N^{\alpha_N} = \varphi (v_1^{\otimes \alpha_1} \ldots v_N^{\otimes \alpha_N}) \, ,
$$
$\alpha_1 + \cdots + \alpha_N = m$. A {\it special basis} of $H^0 (X_K , {\mathcal O} (m))$ is a basis made of such monomials.

\smallskip

Let $r_1 \geq r_2 \geq \cdots \geq r_N$ be $N$ real numbers and ${\bm r} = (r_1 , \ldots , r_N)$. We define the weight of $v_i$ to be $r_i$, the weight of a monomial in $E_K^{\otimes m}$ to be the sum of the weights of the $v_i$'s occuring in it, and the weight of a monomial $u \in H^0 (X_K , {\mathcal O} (m))$ to be the minimum ${wt}_{\bm r} (u)$ of the weights of the monomials in the $v_i$'s mapping to $u$ by $\varphi$. The weight ${wt}_{\bm r} ({\mathcal B})$ of a special basis ${\mathcal B}$ is the sum of the weights of its elements, and $w_{\bm r} (m)$ is the minimum of the weight of a special basis of $H^0 (X_K , {\mathcal O} (m))$.

\smallskip

When $r_1 \geq r_2 \geq \cdots \geq r_N$ are natural integers there exists $e_{\bm r} \in {\mathbb N}$ such that, as $m$ goes to infinity,
$$
w_{\bm r} (m) = e_{\bm r} \, \frac{m^2}{2} + O(m)
$$
(\cite{M}, \cite{Mo} Corollary 3.3).

\smallskip

Our next goal is to find an upper bound for $e_{\bm r}$.

\subsection{}

For every positive integer $i \leq N$ we let $e_i$ be the drop in degree of $X_K$ when projected from ${\mathbb P} (E_K^{\vee})$ to ${\mathbb P} (V_i^{\vee})$. A criterion of Gieseker (\cite{G}, \cite{Mo} Corollary 3.8) tells us that $e_{\bm r} \leq S$ with
$$
S = \underset{1 = i_0 < \ldots < i_{\ell} = N}{\rm min} \ \sum_{j=0}^{\ell - 1} (r_{i_j} - r_{i_{j+1}}) (e_{i_j} + e_{i_{j+1}}) \, .
$$
Note that $S$ is an increasing function in each variable $e_i$. Furthermore, it follows from Clifford's theorem and Riemann-Roch that
\begin{equation}
\label{eq3}
e_i \leq f_i
\end{equation}
for every positive $i \leq N$ -- see \cite{Mo} proof of Theorem~4.4 (N.B.: in \cite{Mo} Theorem~4.4 the filtration of $V_0$ has length $n+1$, while $n = \dim V_0$. In our case, we start the filtration with $V_1$, hence the discrepancy between our definition of $f_i$ and \cite{Mo} {\it loc. cit.}).

\subsection{}

Let $w_1 , \ldots , w_N \in E_K^{\vee}$ be the dual basis of $v_1 , \ldots , v_N$. The linear projection from ${\mathbb P} (V_i^{\vee})$ to ${\mathbb P} (V_{i+1}^{\vee})$ has center the image $\dot w_i$ of $w_i$.

\smallskip

If $y_1 , \ldots , y_N \in E_K^{\vee}$ is the dual basis of $x_1 , \ldots , x_N$, we get
$$
w_i = \left\{ \begin{matrix}
y_i + n_i \, z_i &\mbox{if} \ k \leq i \leq N-3 \\
y_i \hfill &\mbox{else}, \hfill
\end{matrix}
\right.
$$
where $z_i+y_{i+1}$ is a linear combination of $y_{i+2},y_{i+3},\cdots$ with coefficients depending only on
$n_{i+1},n_{i+2},\cdots$.
When $n \ne m$ are two integers, the vectors $y_i + n\,z_i$ and $y_i + m\,z_i$ are linearly independent over $K$, therefore their images in ${\mathbb P} (V_i^{\vee})$ are distinct. Since $e_{N-3} \leq f_{N-3}$ and $g \geq 2$ we get $e_{N-3} \leq d-3$, therefore the image of $X_K$ in ${\mathbb P} (V_i^{\vee})$, $i \leq N-3$, has degree at least 3. Furthermore $\dim {\mathbb P} (V_i^{\vee}) \geq 3$. By Theorem~1 and its Corollary, it follows that we can choose $n_i$ such that $0 \leq n_i < A(d) + d$ and the projection of ${\mathbb P} (V_i^{\vee})$ to ${\mathbb P} (V_{i+1}^{\vee})$ does not change the degree of the image of $X_K$. We fix the integers $n_i$, $k \leq i \leq N-3$, with this property. Hence we have
\begin{equation}
\label{eq4}
e_i = e_k \quad \mbox{whenever} \quad k \leq i \leq N-2 \, .
\end{equation}

\subsection{}

From (\ref{eq3}) and (\ref{eq4}) we conclude that
$$
e_i \leq h_{i,k} \quad \mbox{if} \quad 1 \leq i \leq N
$$
(see 2.3). Hence, by Morrison's main combinatorial theorem, \cite{Mo} Corollary~4.3, for any decreasing sequence of real numbers $r_1 \geq r_2 \geq \cdots \geq r_N$ we have, if $k \geq 2$,
$$
S \leq \psi (\bm r)
$$
with 
$$
\psi (\bm r) = B_k \cdot \sum_{j=1}^N (r_j - r_N) \, .
$$
So, when $r_1 \geq r_2 \geq \cdots \geq r_N = 0$ is a decreasing sequence of real numbers,
$$
e_{\bm r} \leq \psi (\bm r) \, .
$$

From the proof of Theorem~1 in \cite{S} we deduce that, letting
$$
s_i = \log \Vert v_i \Vert - \log \Vert v_N \Vert \, , \quad 1 \leq i \leq N \, ,
$$
\begin{equation}
\label{eq5}
\frac{h(X_K)}{[K:{\mathbb Q}]} + 2d \log \Vert v_N \Vert + \psi (s_1 , s_2 , \ldots , s_{N-1} , 0) \geq 0 \, .
\end{equation}
From (\ref{eq2}) above we get
$$
\log \Vert v_i \Vert \leq \log \Vert x_{i-1} \Vert + \log (1+n_i) \quad \mbox{if} \quad k+1 \leq i \leq N-2
$$
and $\log \Vert v_i \Vert = \log \Vert x_i \Vert$ otherwise.

\smallskip

Since $\log \Vert x_i \Vert = \mu_{N+1-i}$ and $n_i < A+d$ we deduce that

\begin{equation}
\label{eq6}
\psi (s_1 , s_2 , \ldots , s_{N-1} , 0) \leq B_k \biggl( \sum_{i=1}^N (\mu_i - \mu_1) + \mu_{N+1-k} - \mu_3 + (N-2-k) \log (A+d) \biggl) \, .
\end{equation}

From (\ref{eq1}), (\ref{eq5}) and (\ref{eq6}) it follows that
\begin{equation}
\label{eq7}
\frac{h(X_K)}{[K:{\mathbb Q}]} + 2d \, \mu_1 + B_k (N(\mu - \mu_1) + \mu_{N+1-k} - \mu_3) + C \geq 0
\end{equation}
for some constant $C = C(d)$. Since $\mu_3 \geq \mu_1$ the inequality in Theorem~2 follows from (\ref{eq7}).

.

soule@ihes.fr
\end{document}